\documentclass[reqno,11pt]{amsart}

\usepackage{amssymb}
\usepackage[all,2cell,cmtip]{xy}
\UseAllTwocells

\usepackage{comment}
\newtheorem{Proposition}{Proposition}[section]
\newtheorem{Corollary}[Proposition]{Corollary}
\newtheorem{Definition}[Proposition]{Definition}
\newtheorem{Lemma}[Proposition]{Lemma}
\newtheorem{Theorem}[Proposition]{Theorem}
\newtheorem{MainTheorem}{Theorem}

\DeclareMathOperator{\Hom}{Hom}
\DeclareMathOperator{\Sym}{Sym}
\DeclareMathOperator{\Cont}{Cont}
\DeclareMathOperator{\spt}{spt}
\DeclareMathOperator{\spa}{span}
\newcommand{\R}{\mathbb{R}}
\newcommand{\C}{\mathbb{C}}

\newcommand{\ContMan}{\mathbf{ContMan}}
\newcommand{\Vect}{\mathbf{Vect}}
\newcommand{\Nat}{\mathbf{Nat}}
\newcommand{\ed}{\mathrm{d}}
\newcommand{\D}{\mathrm{D}}
\newcommand{\Lo}{\mathrm{L}}
\newcommand{\Po}{\mathrm{P}}
\newcommand{\po}{\mathrm{p}}
\newcommand{\largewedge}{\mbox{\Large $\wedge$}}
\title{Natural operations on differential forms on contact manifolds}
\author{Andreas Bernig}
 
\email{bernig@math.uni-frankfurt.de}

\address{Institut f\"ur Mathematik, Goethe-Universit\"at Frankfurt,
Robert-Mayer-Str. 10, 60054 Frankfurt, Germany}

\thanks{Supported by DFG grants BE 2484/5-1 and BE 2484/5-2.\\ AMS 2010 {\it Mathematics subject
classification}: 
53D10, %Contact manifolds, general
58A10 %Differential forms
}

\begin{document}
\begin{abstract}
We characterize all natural linear operations between spaces of differential forms on contact manifolds. Our main theorem says roughly that such operations are built from some algebraic operators which we introduce and the exterior derivative. 
 \end{abstract}

\maketitle 
%\tableofcontents
%---------------------------------------------------------------------

\section{Introduction}

A classical theorem due to Palais \cite{palais59} characterizes those linear operations on differential forms on a manifold which are compatible with diffeomorphisms. The result is roughly that only the identity and the exterior differential have these properties. The linearity assumption was removed by Kol{\'a}{\v{r}}-Michor-Slov{\'a}k in 1993 \cite{kolar_michor_slovak}. 
A very recent result by Navarro-Sancho \cite{navarro_sancho} generalizes this theorem further by considering natural operations on $k$-tuples of differential forms. They prove that such operations can be written as polynomials in the given forms and their exterior differentials. A special case of this theorem was shown earlier by Freed-Hopkins \cite{freed_hopkins}.

A natural general question in this context is the following. Assume that $M$ is endowed with some extra structure (for instance, a contact or symplectic structure, an almost complex structure etc.). What are the operations on differential forms which are compatible with the diffeomorphisms of $M$ respecting the extra structure? In the present paper, we study this question for contact manifolds and contactomorphisms. 

Before describing more precisely our result, let us recall the definition of contact manifolds, and the construction of the Rumin differential operator. 

A smooth manifold $M$ of dimension $2n+1$, endowed with a completely non-integrable field $Q$ of hyperplanes (called contact planes) is called a contact manifold. We refer to \cite{geiges_book} for contact manifolds and to \cite{cannas_da_silva_book, huybrechts05} for some of the symplectic linear algebra which we will use. 

Locally, the contact plane can be written as $\ker \alpha$ for some $1$-form $\alpha$. The condition of complete non-integrability is that $\alpha \wedge (\ed \alpha)^n \neq 0$. A contactomorphism (resp. local contactomorphism) $\phi:(M,Q) \to (\tilde M,\tilde Q)$ is a diffeomorphism (resp. local diffeomorphism) such that $\ed\phi(Q)=\tilde Q$. If $Q=\ker \alpha, \tilde Q=\ker \tilde \alpha$, then a (local) diffeomorphism $\phi$ is a (local) contactomorphism if and only if $\phi^*\tilde \alpha=f \alpha$ for some nowhere vanishing function $f$. If $f \equiv 1$, $\phi$ is called a strict contactomorphism. The group of all contactomorphisms of $(M,Q)$ is denoted by $\Cont(M)$.

The standard example of a contact manifold is $\R^{2n+1}$ with coordinates $(x_1,\ldots,x_n,y_1,\ldots,y_n,z)$ and $\alpha:=\ed z +\sum_{i=1}^n x_i \ed y_i$. Darboux' theorem states that locally each contact manifold of dimension $2n+1$ is contactomorph to $\R^{2n+1}$ with its standard contact structure.

Let $\Omega^a(M)$ denote the space of differential $a$-forms on $M$. Rumin constructed an operator $\D:\Omega^n(M) \to \Omega^{n+1}(M)$ as follows. Let $\omega \in \Omega^n(M)$. Restrict $\ed\omega$ to the contact plane. By some basic linear symplectic algebra, we may write $\ed\omega|_Q=- d\alpha \wedge \xi|_Q$, where $\xi \in \Omega^{n-1}(M)$. Moreover, $\alpha \wedge \xi$ is unique and independent of the choice of the local defining $1$-form $\alpha$. Then 
\begin{displaymath}
 \D(\omega):=\ed(\omega+\alpha \wedge \xi)
\end{displaymath}
is well-defined. The second order differential operator $\D$ is called Rumin operator. Another natural operator is given by $\mathrm{Q}\omega:=\omega+\alpha \wedge \xi$, so that $\D=\ed \circ \mathrm{Q}$. By construction, both $\D$ and $\mathrm{Q}$ are equivariant with respect to contactomorphisms, i.e. $\D(\phi^*\omega)=\phi^*\D\omega$ and similarly for $\mathrm{Q}$. The Rumin operator fits into a sequence of differential operators defining the so-called Rumin-de Rham complex, whose cohomology is isomorphic to the usual de Rham cohomology \cite{rumin94}. For applications of the Rumin operator in the theory of valuations, see the notes below. 

Let us now construct some more operators, which may be seen as refinements of Rumin's construction and which will play an essential role in the following. We use a locally defining $1$-form $\alpha$, and will show later on that the construction is independent of this choice.

Let 
\begin{displaymath}
 \Lo:\Omega^*(M) \to \Omega^{*+2}(M), \omega \mapsto \omega \wedge \ed\alpha
\end{displaymath}
be the Lefschetz operator. 

A form $\pi \in \Omega^i(M)$ is called primitive if $i \leq n$ and $\Lo^{n-i+1}\pi|_Q=0$. Any form $\omega \in \Omega^a(M)$ can be decomposed as 
\begin{equation}
 \omega|_Q=\sum_{\substack{0 \leq i \leq \min(a,2n-a)\\i \equiv a \mod 2}} \Lo^{\frac{a-i}{2}} \pi_{i}|_Q
\end{equation}
with $\pi_i \in \Omega^i(M)$ primitive, and $\pi_i|_Q$ is unique. 

We define for $0 \leq i \leq \min(a-2,2n-a), i \equiv a \mod 2$ the maps 
\begin{displaymath}
 \Po_{a,i}:\Omega^a(M) \to \Omega^{a-1}(M), \omega \mapsto \Lo^{\frac{a-i-2}{2}} \pi_{i} \wedge \alpha.
\end{displaymath}

Then $\Po_{a,i}$ is compatible with contactomorphisms. The Rumin operator may be written as 
\begin{displaymath}
\D=\ed+(-1)^n\sum_{\substack{0 \leq i \leq n-1\\i \equiv n-1 \mod 2}} \ed \circ \Po_{n+1,i} \circ \ed.
\end{displaymath}

Let us now state our main theorem. Denote by $\ContMan_{2n+1}$ the category whose objects are $(2n+1)$-dimensional contact manifolds and whose morphisms are local contactomorphisms. Let $\mathbf{Vect}$ be the category whose objects are real vector spaces and whose morphisms are linear maps. 

Let $\Omega^a:\ContMan_{2n+1} \to \mathbf{Vect}$ be the contravariant functor  which assigns to each contact manifold $(M,Q)$ the vector space $\Omega^a(M)$ of $a$-forms on $M$, and to each local contactomorphism $\phi:(M,Q) \to (M',Q')$ the linear map $\phi^*:\Omega^a(M') \to \Omega^a(M)$.

\begin{Definition}
A natural linear contact morphism $\Po:\Omega^a \to \Omega^b$ is a natural transformation of functors  
\begin{displaymath}
\xymatrix@C+2pc{
\ContMan_{2n+1} \rtwocell<5>^{\Omega^a}_{\Omega^b}{P} & \Vect
}
\end{displaymath}
The vector space of all natural linear contact morphisms $\Po:\Omega^a \to \Omega^b$ is denoted by $\Nat^{2n+1}_{a,b}$. 
\end{Definition}

This means that for each contact manifold $(M,Q)$ of dimension $2n+1$, we are given a linear map $\Po:\Omega^a(M) \to \Omega^b(M)$; and if $\phi:M \to \tilde M$ is a local contactomorphism, then $\Po(\phi^*\omega)=\phi^* \Po(\omega)$ for all $\omega \in \Omega^a(\tilde M)$.

\begin{MainTheorem} \label{mainthm_fixed_dimension}
 The space $\Nat^{2n+1}_{a,b}$ of natural linear contact morphisms $\Po:\Omega^a \to \Omega^b$ is trivial if $b \neq a-1,a,a+1$ and is given by 
 \begin{displaymath}
\Nat^{2n+1}_{a,b}= \begin{cases}
  \spa \{\Po_{a,i}\} & \text{ if } b=a-1,\\
  \spa\{\mathrm{d} \circ \Po_{a,i}, \Po_{a+1,i} \circ \ed, \mathrm{id}\} & \text{ if } b=a,\\
  \spa\{\ed, \ed \circ \Po_{a+1,i} \circ \ed\} & \text{ if } b=a+1.
 \end{cases}
\end{displaymath}
 Here $i$ ranges over those indices where $\Po_{*,i}$ is defined. 
\end{MainTheorem}

The first step of the proof is to apply Peetre's theorem to show that $\Po$ must be a linear differential operator. The description of all equivariant linear differential operators is an algebraic problem, which we solve using the first fundamental theorem of invariant theory for the symplectic group. 

Let us comment on some possible variations of this theorem. 
\begin{enumerate}

\item Palais' theorem characterizes the space of all linear operators $\Po:\Omega^a(M) \to \Omega^b(M)$ which commute with diffeomorphisms on a fixed manifold $M$. If $b>0$, then this space is spanned by $\mathrm{id}$ if $b=a$ and by $\ed$ if $b=a+1$. If $b=0$ and assuming that $M$ is compact, oriented and does not admit an orientation reversing diffeomorphism, the map $\omega \mapsto \int_M \omega \in \R \subset \Omega^0(M)$ is another such operator. 

In a similar spirit, we could study linear operators $\Po:\Omega^a(M) \to \Omega^b(M)$ which commute with contactomorphisms on a fixed contact manifold $(M,Q)$. If $b>0$, then one can use Darboux' theorem and homotheties of $\R^{2n+1}$ to prove that such operators are support non-increasing. Then our proof of Theorem \ref{mainthm_fixed_dimension} yields that only the operators from Theorem \ref{mainthm_fixed_dimension} appear. For $b=0$ there may be non-local operators but we do not have a classification of them.  

\item Instead of local contactomorphisms, we could use strict local contactomorphisms, i.e. diffeomorphisms which leave a contact form $\alpha$ invariant. In this case, there are other natural linear operators, like $\omega \mapsto \alpha \wedge \ed \alpha \wedge \ed \omega$ or $\omega \mapsto \mathcal{L}_T\omega$, where $T$ is the Reeb vector field associated to $\alpha$. We do not have a complete description in this case. 
 \item  Let $\ContMan$ be the category of contact manifolds and contact embeddings. We denote by $\Omega^a:\ContMan \to \Vect$ the contravariant functor which associates to a contact manifold $(M,Q)$ the space $\Omega^a(M)$ and to each contact embedding $\iota:(M,Q) \to (\tilde M,\tilde Q)$ the pull-back $\iota^*:\Omega^a(\tilde M) \to \Omega^a(M)$. Here $M$ and $\tilde M$ are not necessarily of the same dimension. Using Theorem \ref{mainthm_fixed_dimension}, one can show that the space of natural transformations $\Po:\Omega^a \to \Omega^b$ over $\ContMan$ is spanned by $\mathrm{id}$ if $a=b$, by the exterior differential $\ed$ if $b=a+1$, and trivial otherwise. In other words, only those operators from Theorem \ref{mainthm_fixed_dimension} are dimension-independent which are natural with respect to the larger collection of smooth maps.   
 \item Another possibility is to drop the linearity assumption. Again, more operators appear, for instance $\omega \mapsto \omega \wedge \omega$ or $\omega \mapsto \omega \wedge \Po_{a,i}\omega \wedge \ed \omega$. Using the Peetre-Slov{\'a}k theorem as in \cite{navarro_sancho}, one can show that a natural operator (satisfying an additional regularity assumption) on $\ContMan_{2n+1}$ must be a differential operator. We do not have a complete description of these operators.
 \item One could also, as in \cite{navarro_sancho}, consider natural operations on $k$-tuples of forms. There are many of them, for instance $(\omega_1,\omega_2) \mapsto \Po_{a_1+a_2,i}(\omega_1 \wedge \omega_2)$ or $(\omega_1,\omega_2) \mapsto \ed(\omega_1 \wedge \Po_{a_2,i}\omega_2)$. 
\end{enumerate}
  
Let us finish this introduction with some remarks on applications of the Rumin operator and of the operator $\mathrm{Q}$ in the theory of smooth valuations on manifolds, which stimulated our interest in this operator and in the problem of its characterization. Alesker has defined a far-reaching theory of valuations on manifolds (see \cite{alesker_survey07} and the references there). Among many other applications, this theory was of fundamental importance in the recent study of integral geometry of complex projective and complex hyperbolic spaces \cite{alesker_bernig, bernig_fu_hig, bernig_fu_solanes, wannerer_area_measures, wannerer_unitary_module}. Given an oriented manifold $X$ of dimension $n$, the cosphere bundle $\mathbb{P}_+(X)$ is a contact manifold of dimension $2n-1$. Given a pair of forms $\omega \in \Omega^{n-1}(\mathbb{P}_+(X)), \phi \in \Omega^n(X)$, one can construct a smooth valuation on $X$ by integration over the conormal cycle. The pairs which induce the trivial valuation were characterized in \cite{bernig_broecker07} as those for which $\D\omega+\pi^*\phi=0, \pi_*\omega=0$. Here $\D$ is the Rumin operator, $\pi: \mathbb{P}_+(X) \to X$ is the natural projection map and $\pi_*$ denotes push-forward. 

The Rumin operator was also a key tool in the construction of a convolution product on smooth and generalized translation invariant valuations on affine space \cite{bernig_faifman, bernig_fu06}, the product of smooth valuations on manifolds \cite{alesker_bernig, fu_alesker_product} (note that the operator $\mathrm{Q}$ played some role in this case), and the operations of push-forward and pull-back of smooth valuations \cite{alesker_intgeo}.

Furthermore, the Rumin-de Rham complex appears naturally in the study of invariant differential operators on parabolic geometries \cite{bryant_eastwood_gover_neusser}. More precisely, the symplectic group $G=\mathrm{Sp}_{2n+2}\R$ acts on the unit sphere $S^{2n+1}$, and the stabilizer is a parabolic subgroup $P$. The Bernstein-Gelfand-Gelfand sequence \cite{calderbank_diemer, cap_slovak_soucek} corresponding to the trivial representation of $G$ is then the Rumin-de Rham complex.

\subsection{Plan of the paper} 
We first construct the operators $\Po_{a,i}$ and show that they indeed commute with local contactomorphisms. We also show how the Rumin operator can be expressed using these operators.  In Section \ref{sec_peetre} we state Peetre's theorem and apply it to our situation. In Section \ref{sec_isotropy_group} we study the infinitesimal action of the group of contactomorphisms on a fixed tangent space. In particular we will see that this action contains the non-compact symplectic group.

In Section \ref{sec_principal_symbols} we study principal symbols of invariant differential operators and finish the proof of the main theorem modulo a classification of such symbols. This classification is contained in Sections \ref{sec_classification_horizontal}, \ref{sec_reduction_horizontal} and \ref{sec_operators_degree2}. The basic idea is to apply the first fundamental theorem for the non-compact symplectic group to the horizontal part of a symbol. We then show that an invariant symbol is already characterized by its restriction to the horizontal part. Not all symbols which are invariant under the symplectic group will be invariant under the full group of infinitesimal contactomorphisms, and in Section \ref{sec_operators_degree2} we prove that only the symbols of the operators from Theorem \ref{mainthm_fixed_dimension} have this property.

\subsection*{Acknowledgments}
The problem of a characterization of the Rumin operator in the spirit of Palais' theorem was brought to my attention by Marc Burger back in 2005. In the course of working out this material, I profited from discussions with Semyon Alesker, Dmitry Faifman and Thomas Mettler. I thank Myhailo Saienko for many useful remarks on a first draft of this manuscript.

\section{Construction of natural contact morphisms}

Let $(M,Q)$ be a contact manifold of dimension $2n+1$. In the introduction, we have defined operators $\Po_{a,i}:\Omega^a(M) \to \Omega^{a-1}(M)$, where $0 \leq i \leq \min(a-2,2n-a), i \equiv a \mod 2$. 

\begin{Lemma} \label{lemma_naturality}
 The operator $\Po_{a,i}$ is linear and equivariant under local contactomorphisms, i.e. an element of $\mathbf{Nat}^{2n+1}_{a,a-1}$.  
\end{Lemma}

\proof
It is clear that $\Po_{a,i}$ is linear and equivariant under strict local contactomorphisms. It therefore remains to show that it is independent of the choice of $\alpha$. 
 
Let $\tilde \alpha=f \alpha$ be another locally defining $1$-form, where $f$ is a nowhere vanishing smooth function. Let $\tilde \Lo$ denote the corresponding Lefschetz operator and $\tilde \Po_{a,i}$ the operator corresponding to this choice of a contact form. Then 
\begin{displaymath}
 \omega|_Q=\sum_{\substack{0 \leq i \leq \min(a,2n-a)\\ i \equiv a \mod 2}} \tilde \Lo^{\frac{a-i}{2}} \tilde \pi_{i}|_Q,
\end{displaymath}
with $\tilde \pi_i:=f^{-\frac{a-i}{2}} \pi_i$ primitive. Hence 
\begin{align*}
\tilde \Po_{a,i}(\omega) & = \tilde \Lo^{\frac{a-i-2}{2}} \tilde \pi_{i} \wedge \tilde \alpha\\
& = (f \ed \alpha + \ed f \wedge \alpha)^\frac{a-i-2}{2} \wedge f^{-\frac{a-i}{2}} \pi_i \wedge f\alpha\\
& = \ed \alpha^\frac{a-i-2}{2} \wedge \pi_i \wedge \alpha\\
& = \Lo^{\frac{a-i-2}{2}} \pi_{i} \wedge \alpha\\
& = \Po_{a,i}(\omega).
\end{align*}
\endproof

Composing with the exterior differential $\ed$ yields additional natural operators. However, the next lemma shows that the operators given in Theorem \ref{mainthm_fixed_dimension} are the only relevant ones. 

\begin{Lemma}
 For all $a,i,j$ we have 
 \begin{displaymath}
  \Po_{a-1,j} \circ \Po_{a,i} = 0
 \end{displaymath}
 and 
 \begin{displaymath}
  \Po_{a,j} \circ \ed \circ \Po_{a,i}=
  \begin{cases} 0 & i \neq j\\
  (-1)^a \Po_{k,i} & i=j.
  \end{cases}
\end{displaymath}
\end{Lemma}

\proof
The first equation is immediate from the fact that $\Po_{a,i} \omega|_Q=0$ for all $\omega$.

For the second equation, let 
\begin{displaymath}
 \omega|_Q=\sum_{\substack{0 \leq i \leq \min(a,2n-a)\\ i \equiv a \mod 2}} \Lo^{\frac{a-i}{2}} \pi_i
\end{displaymath}
be the Lefschetz decomposition of $\omega \in \Omega^a(M)$. Then 
\begin{displaymath}
 \ed \circ \Po_{a,i} \omega=\ed(\Lo^{\frac{a-i-2}{2}} \pi_i \wedge \alpha) \equiv (-1)^i \Lo^{\frac{a-i}{2}} \pi_i \mod \alpha.
\end{displaymath}
The Lefschetz decomposition of $(\ed \circ \Po_{a,i} \omega)|_Q$ is thus given by $(-1)^i \Lo^{\frac{a-i}{2}} \pi_i$. It follows that $\Po_{a,j} \circ \ed \circ \Po_{a,i} \omega=0$ if $i \neq j$. For $i=j$ we obtain 
\begin{displaymath}
 \Po_{a,i} \circ \ed \circ \Po_{a,i} \omega=\Po_{a,i}((-1)^i\Lo^{\frac{a-i}{2}} \pi_i)=(-1)^i\Lo^{\frac{a-i-2}{2}} \pi_i \wedge \alpha=(-1)^i\Po_{a,i}\omega.
\end{displaymath}
\endproof

\begin{Lemma} \label{lemma_compositions}
The operator $\mathrm{Q}$ is given by 
\begin{displaymath}
 \mathrm{Q}=\mathrm{id}+(-1)^n \sum_{\substack{0 \leq i \leq n-1\\i \equiv n-1 \mod 2}} \Po_{n+1,i} \circ \ed.
\end{displaymath}
The Rumin differential $\D:\Omega^n(M) \to \Omega^{n+1}(M)$ is given by
\begin{displaymath}
\D=\ed+(-1)^n \sum_{\substack{0 \leq i \leq n-1\\i \equiv n-1 \mod 2}} \ed \circ \Po_{n+1,i} \circ \ed.
\end{displaymath}
\end{Lemma}

\proof
Let $\omega \in \Omega^n(M)$. Let 
\begin{displaymath}
\ed \omega|_Q=\sum_{\substack{0 \leq i \leq n-1\\i \equiv n-1 \mod 2}} \Lo^{\frac{n+1-i}{2}} \pi_i
\end{displaymath}
be the Lefschetz decomposition of $\ed\omega$.
It follows that 
\begin{displaymath}
\ed \omega|_Q=- \Lo \xi|_Q
\end{displaymath}
with 
\begin{displaymath}
 \xi:=- \sum_{\substack{0 \leq i \leq n-1\\i \equiv n-1 \mod 2}} \Lo^{\frac{n-1-i}{2}} \pi_{i},
\end{displaymath}
which implies that $\mathrm{Q}\omega=\omega+\alpha \wedge \xi$ and $\D\omega=\ed(\omega+\alpha \wedge \xi)$. The statement follows from 
\begin{align*}
 \alpha \wedge \xi & =(-1)^{n-1} \xi \wedge \alpha\\
 & = (-1)^n \sum_{\substack{0 \leq i \leq n-1\\i \equiv n-1 \mod 2}} \Lo^{\frac{n-i-1}{2}} \pi_{i} \wedge \alpha\\
 & = (-1)^n \sum_{\substack{0 \leq i \leq n-1\\i \equiv n-1 \mod 2}} \Po_{n+1,i} \circ \ed \omega.
\end{align*}
\endproof

\section{Peetre's theorem}
\label{sec_peetre}

\begin{Definition}
 Let $E,F$ be smooth vector bundles over $M$ and denote by $\mathcal{E}(M),\mathcal{F}(M)$ the spaces of smooth global sections. A linear operator $\Po:\mathcal{E}(M) \to \mathcal{F}(M)$ is called support non-increasing if $\mathrm{spt}(\Po s) \subset \mathrm{spt}(s)$ for each $s \in \mathcal{E}(M)$. 
\end{Definition}

\begin{Theorem}[Peetre \cite{peetre}]  \label{thm_peetre}
Let $E,F$ be smooth vector bundles over $M$ and let $\Po:\mathcal{E}(M) \to \mathcal{F}(M)$ be a linear, support non-increasing operator. Then, for every compact subset $K \subset M$, $\Po$ is a differential operator on $K$ of finite order.
\end{Theorem}

We refer to \cite{kolar_michor_slovak} for more information on this theorem and its generalizations. 

\begin{Proposition} \label{prop_finite_order}
Let $\R^{2n+1}$ be endowed with its standard contact structure. Let $\Po:\Omega^a \to \Omega^b$ be a natural linear contact morphism. Then $\Po|_{\R^{2n+1}}$ is a differential operator of finite order. 
\end{Proposition}

\proof
We are going to apply Peetre's theorem \ref{thm_peetre}. For this we have to check that $\Po$ is support non-increasing. 

Let $\omega \in \Omega^a(\R^{2n+1})$ and $p \notin \spt \omega$. Let $U$ be an open neighborhood of $p$ disjoint from $\spt \omega$. We may regard $U$ as a contact manifold with the induced contact structure. Let $\iota:U \to \R^{2n+1}$ be the inclusion. By naturality, 
\begin{displaymath}
\iota^* \Po \omega=\Po \iota^*\omega=0, 
\end{displaymath}
hence $\spt \Po \omega \cap U = \emptyset$. It follows that $\Po$ is support non increasing. By Peetre's theorem, $\Po$ is a linear differential operator of finite degree on each compact subset $K$ of $\R^{2n+1}$. Since the group of contactomorphisms acts transitively on $\R^{2n+1}$, the order is finite on all of $\R^{2n+1}$ and we are done. 
\endproof

\section{The action of the isotropy group}
\label{sec_isotropy_group}

Consider $\R^{2n+1}$ with its standard contact form
\begin{displaymath}
 \alpha:=\ed z+\sum_{i=1}^n x_i \ed y_i
\end{displaymath}
and let $V:=T_0\R^{2n+1} \cong \R^{2n+1}$ be the tangent space at $0$.

The associated Reeb vector field is defined by the two equations $i_R \alpha \equiv 1$ and $i_R \ed \alpha \equiv 0$, which implies that $R=\frac{\partial}{\partial z}$.

Let $W:=Q_0$ be the span of $\left.\frac{\partial}{\partial x_i}\right|_0,\left.\frac{\partial}{\partial y_i}\right|_0, i=1,\ldots,n$. The restriction of $\ed \alpha|_0$ to $W$ defines a symplectic form $\Omega$ on $W$. Let $D$ be the span of $\left.\frac{\partial}{\partial z}\right|_0$. Then $V \cong W \oplus D$. 

Infinitesimal contactomorphisms are called contact isotopies. They are induced by vector fields $X$ with  
\begin{displaymath}
 \alpha(X)=h, \quad i_X \ed\alpha=\ed h(R)\alpha-\ed h,
\end{displaymath}
where $h$ is a smooth function on $\R^{2n+1}$ \cite[Theorem 2.3.1]{geiges_book}. Solving these equations gives 
\begin{equation} \label{eq_contact_field}
 X=\sum_{i=1}^n \left(\frac{\partial h}{\partial z}x_i-\frac{\partial h}{\partial y_i}\right) \frac{\partial}{\partial x_i}+\sum_{i=1}^n \frac{\partial h}{\partial x_i} \frac{\partial}{\partial y_i}+\left(h-\sum_{i=1}^n \frac{\partial h}{\partial x_i}x_i\right) \frac{\partial}{\partial z}.
\end{equation}
Such vector fields are called contact vector fields.

In the following, we identify $\R^{2n+1}=\C^n \times \R, (x,y,z) \mapsto (x+iy,z)$. Let $\C^n$ be endowed with its standard symplectic form and let $\mathrm{Sp}_{2n}\R$ be the Lie group of linear automorphisms of $\C^n$ preserving the symplectic form. 

\begin{Proposition} \label{prop_stabilizer}
Let $G$ be the stabilizer of $\Cont(\R^{2n+1})$ at $0$, which acts naturally on  the tangent space $V=T_0\R^{2n+1}=W \oplus D \cong \C^n \oplus \R$. Then $G$ contains the following three subgroups: 
\begin{enumerate}
 \item the (non-compact) symplectic group $G_0:=\mathrm{Sp}_{2n}\R$ acting on $\R^{2n+1}=\C^n \oplus \R$ in the standard way on the first factor and trivially on the second factor. 
 \item The group $\R_{>0}$, acting by $(x+iy,z) \mapsto (\lambda (x+iy),\lambda^2 z), \lambda>0, (x+iy,z) \in \C^n \oplus \R$.
 \item The group $\R^{2n}$, acting by $(x+iy,z) \mapsto (x+zv_1+i(y+zv_2),z), (v_1,v_2) \in \R^{2n}, (x+iy,z) \in \C^n \oplus \R$.
\end{enumerate}
\end{Proposition}

\proof
Since the three displayed groups are connected, it is enough to prove the statement on the level of Lie algebras. 

Let $\phi_t$ be the flow generated by a contact vector field $X$. Let $h \in C^\infty(\R^{2n+1})$ be the corresponding function from \eqref{eq_contact_field}. Suppose that $\phi_t(0)=0$ for all $t$. This means that $X|_0=0$, which is equivalent to 
\begin{equation} \label{eq_conditions_h}
\left. \frac{\partial h}{\partial x_i}\right|_0=0,\left.\frac{\partial h}{\partial y_i}\right|_0=0, h(0)=0. 
\end{equation}

Writing $\phi=\left(\begin{array}{c} \phi^1\\ \vdots \\ \phi^{2n+1}\end{array}\right)$, we obtain 
\begin{displaymath}
 d(\phi_t)=\left(\begin{array}{c c c c c c c}  \frac{\partial \phi^1_t}{\partial x_1} & \ldots & \frac{\partial \phi^1_t}{\partial x_n} & \frac{\partial \phi^1_t}{\partial y_1} & \ldots & \frac{\partial \phi^1_t}{\partial y_n} & \frac{\partial \phi^1_t}{\partial z}\\
 \frac{\partial \phi_t^2}{\partial x_1} & \ldots & \frac{\partial \phi_t^2}{\partial x_n} & \frac{\partial \phi_t^2}{\partial y_1} & \ldots & \frac{\partial \phi_t^2}{\partial y_n} & \frac{\partial \phi_t^2}{\partial z}\\
 \vdots & \vdots & \vdots & \vdots & \vdots & \vdots & \vdots\\
 \frac{\partial \phi_t^{2n+1}}{\partial x_1} & \ldots & \frac{\partial \phi_t^{2n+1}}{\partial x_n} & \frac{\partial \phi_t^{2n+1}}{\partial y_1} & \ldots & \frac{\partial \phi_t^{2n+1}}{\partial y_n} & \frac{\partial \phi_t^{2n+1}}{\partial z}              
\end{array}\right)
\end{displaymath}

Define $n \times n$-matrices $A,B,C$ by 
\begin{displaymath}
 A := \left(\left.\frac{\partial^2 h}{\partial x_i \partial x_j}\right|_0\right)_{i,j}, B := \left(\left.\frac{\partial^2 h}{\partial x_i \partial y_j}\right|_0\right)_{i,j}, C := \left(\left.\frac{\partial^2 h}{\partial y_i \partial y_j}\right|_0\right)_{i,j}
\end{displaymath}
and $n$-vectors
\begin{displaymath}
 E_x  := \left(\left.\frac{\partial^2 h}{\partial x_i \partial z}\right|_0\right)_{i}, E_y := \left(\left.\frac{\partial^2 h}{\partial y_i \partial z}\right|_0\right)_{i}.
\end{displaymath}

Taking the derivative at $t=0$, changing the order of derivatives and evaluating at $0$ gives us 
\begin{displaymath}
\left.\frac{d}{dt}\right|_{t=0} d(\phi_t)=\left(\begin{array}{c c c}  
\frac{\partial h}{\partial z} \mathrm{Id}_n -B^t & -C & -E_y \\
A & B & E_x \\
0 & 0 & \left.\frac{\partial h}{\partial z}\right|_0          
\end{array}\right)
\end{displaymath}

Since $h$ is an arbitrary function satisfying \eqref{eq_conditions_h}, the matrices $A,B,C$, the vectors $E_x,E_y$ and the real number $\left.\frac{\partial h}{\partial z}\right|_0$ are arbitrary as well, with the obvious restriction that  $A$ and $C$ are symmetric. 

To obtain the Lie algebras of the subgroups displayed in the proposition, we choose $E_x=E_y=0, \left.\frac{\partial h}{\partial z}\right|_0=0$ for $G_0$; $B=\mathrm{Id}, \left.\frac{\partial h}{\partial z}\right|_0=2, E_x=E_y=0$ for $\R_{>0}$; and $A=B=C=0, \left.\frac{\partial h}{\partial z}\right|_0=0$ for $\R^{2n}$.
This proves the proposition.
\endproof

\section{Principal symbols}
\label{sec_principal_symbols}

Let $\Po:\Omega^a \to \Omega^b$ be a natural linear contact morphism. By Proposition \ref{prop_finite_order}, $\Po|_{\R^{2n+1}}$ is a linear differential operator of finite order $r$. Let 
\begin{displaymath}
 \sigma_{\Po}:\largewedge^a T^*\R^{2n+1} \otimes \Sym^r(T^*\R^{2n+1}) \to \largewedge^b T^*\R^{2n+1}  
\end{displaymath}
be the principal symbol of $\Po|_{\R^{2n+1}}$. Then $\sigma_{\Po}$ is $\Cont(\R^{2n+1})$-equivariant and
\begin{displaymath}
 \sigma_{\Po,0}:\largewedge^a T^*_0\R^{2n+1} \otimes \Sym^r(T^*_0\R^{2n+1}) \to \largewedge^b T^*_0\R^{2n+1}
\end{displaymath}
is $G$-equivariant, where $G$ is the stabilizer of $\Cont(\R^{2n+1})$ at the point $0$.  

Let us describe the symbols of the operators from Theorem \ref{mainthm_fixed_dimension}. Set $V:=T_0\R^{2n+1}$ and let $\alpha_0 \in V^*$ and $\Omega_0 \in \largewedge^2 V^*$ be the values of the contact form $\alpha$ and of the form $\Omega=\ed\alpha$ at $0$. Let $\Lo:\largewedge^*V^* \to \largewedge^{*+2}V^*$ be multiplication by $\Omega_0$. 

We call an element $\pi \in \largewedge^i V^*$ primitive if $i \leq n$ and $\Lo^{n-i+1}\pi|_{Q_0}=0$. 

If $\phi \in \largewedge^a V^*$, then the restriction to $Q_0$ can be decomposed as 
\begin{equation}
 \phi|_{Q_0}=\sum_{\substack{0 \leq i \leq \min(a,2n-a)\\i \equiv a \mod 2}} \Lo^{\frac{a-i}{2}} \pi_{i}
\end{equation}
with $\pi_i \in \largewedge^i V^*$ primitive, and $\pi_i|_{Q_0} \in \largewedge^i Q_0^*$ is unique. Define operators $\po_{a,i}: \largewedge^a V^* \to \largewedge^{a-1} V^*$ for $0 \leq i \leq \min(a-2,2n-a), i \equiv a \mod 2$ by 
 \begin{displaymath}
  \po_{a,i}(\phi):=\Lo^{\frac{a-i-2}{2}} \pi_i \wedge \alpha_0.
 \end{displaymath}

We then obtain the following table of symbols 
\begin{displaymath}
\begin{array}{c | c}
\textbf{ operator} \Po & \textbf{symbol } \sigma_{P,p}  \\ \hline
\Po_{a,i} & \po_{a,i} \\
\ed \circ \Po_{a,i} & \phi \otimes \tau \mapsto   \po_{a,i}(\phi) \wedge \tau\\
\Po_{a+1,i} \circ \ed & \phi \otimes \tau \mapsto \po_{a+1,i}(\phi \wedge \tau)\\
\mathrm{id} & \mathrm{id}\\
\ed & \phi \otimes \tau \mapsto \phi \wedge \tau\\
\ed \circ \Po_{a+1,i} \circ \ed & \phi \otimes \tau_1 \otimes \tau_2 \mapsto \po_{a+1,i}(\phi \wedge \tau_1) \wedge \tau_2+\po_{a+1,i}(\phi \wedge \tau_2) \wedge \tau_1  
\end{array}
\end{displaymath}

Define the vector space 
\begin{displaymath}
R_{a,b}^{2n+1}:=\begin{cases}
  \spa\{\Po_{a,i}\} & \text{ if } b=a-1\\
  \spa\{\ed \circ \Po_{a,i}, \Po_{a+1,i} \circ \ed, \mathrm{id}\} & \text{ if } b=a\\
  \spa\{\ed, \ed \circ \Po_{a+1,i} \circ \ed\} & \text{ if } b=a+1\\
  0 & \text{ if } b \neq a-1,a,a+1.
 \end{cases}
\end{displaymath}
Note that the operators on the right hand side are those appearing in Theorem \ref{mainthm_fixed_dimension}. The theorem can thus be restated as $\mathbf{Nat}^{2n+1}_{a,b}=R_{a,b}^{2n+1}$ for all $a,b$. The inclusion $\supset$ follows from Lemma \ref{lemma_naturality}.

\begin{Proposition} \label{prop_invariant_symbols}
Let $\R^{2n+1}$ be endowed with its standard contact structure and let $G$ be the stabilizer of the point $0 \in \R^{2n+1}$. Let $\sigma:\largewedge^a T^*_0\R^{2n+1} \otimes \Sym^r(T^*_0\R^{2n+1}) \to \largewedge^b T^*_0\R^{2n+1}$ be $G$-equivariant. Then there exists $\Po \in R^{2n+1}_{a,b}$ such that $\sigma_{\Po,0}=\sigma$. 
\end{Proposition}

We postpone the proof of this proposition to the following sections and show how Theorem \ref{mainthm_fixed_dimension} follows from it. 

\proof[Proof of Theorem \ref{mainthm_fixed_dimension}]
We use induction on $r$ to show that each $\Cont(\R^{2n+1})$-equivariant linear differential operator $\Po:\Omega^a(\R^{2n+1}) \to \Omega^b(\R^{2n+1})$ of order $\leq r$ equals the restriction to $\R^{2n+1}$ of an element from $R^{2n+1}_{a,b}$. The induction start is the empty case $r=-1$ (with the convention that a differential operator of degree $\leq -1$ is trivial). 

For the induction step, let $r \geq 0$ and let $\Po$ be a linear equivariant differential operator of order $\leq r$. Then $\sigma_{\Po,0}$ is $G$-equivariant. By Proposition \ref{prop_invariant_symbols} there exists an operator $\Po_1 \in R^{2n+1}_{a,b}$ such that $\sigma_{\Po,0}=\sigma_{\Po_1,0}$. Since both $\Po$ and $\Po_1$ are $\Cont(\R^{2n+1})$-equivariant and $\Cont(\R^{2n+1})$ acts transitively on $\R^{2n+1}$, we even have $\sigma_{\Po,q}=\sigma_{\Po_1,q}$ for all $q \in \R^{2n+1}$. It follows that $\Po-\Po_1$ is a $\Cont(\R^{2n+1})$-equivariant linear differential operator of degree $\leq r-1$. By induction hypothesis, $\Po-\Po_1 \in R^{2n+1}_{a,b}$ and we are done. 

Let us now finish the proof of the theorem. Let $\Po:\Omega^a \to \Omega^b$ be a natural linear contact morphism. By Proposition \ref{prop_finite_order}, $\Po|_{\R^{2n+1}}$ is a $\Cont(\R^{2n+1})$-equivariant linear differential operator of finite order. By what we have shown, there exists $\Po_1 \in R^{2n+1}_{a,b}$ with $\Po|_{\R^{2n+1}}=\Po_1|_{\R^{2n+1}}$. By naturality, if $V \subset \R^{2n+1}$ is an open subset, considered as a contact submanifold, then $\Po|_V=\Po_1|_V$.

We claim that $\Po|_M=\Po_1|_M$ for every contact manifold $M$ of dimension $2n+1$. Let $p \in M$. By Darboux' theorem, there exists an open neighborhood $U$ of $p$ in $M$, an open neighborhood $V$ of $0$ in $\R^{2n+1}$ and a contactomorphism $\phi:V \to U$. Since $\Po$ and $\Po_1$ are natural with respect to contactomorphisms, we have 
\begin{displaymath}
 \phi^*\Po \omega=\Po \phi^*\omega=\Po_1 \phi^* \omega=\phi^* \Po_1\omega.
\end{displaymath}
Hence $\Po \omega$ and $\Po_1 \omega$ agree in a neighborhood of $p$. Since $p$ was arbitrary, we conclude that $\Po \omega=\Po_1 \omega$ on $M$, which finishes the proof.
\endproof

\section{Classification of invariant horizontal symbols}
\label{sec_classification_horizontal}

Let $(W,\Omega)$ be a symplectic vector space of dimension $2n$. Let $G_0:=\mathrm{Sp}(W,\Omega)=\mathrm{Sp}_{2n}\R$ be the symplectic group. 

Let $\Lo:\largewedge^* W^* \to \largewedge^{*+2}W^*, \omega \mapsto \omega \wedge \Omega$ be the Lefschetz operator. As before, an element $\phi \in \largewedge^aW^*$ can be uniquely decomposed as 
\begin{equation} \label{eq_lefschetz_decomp}
 \phi=\sum_{\substack{0 \leq i \leq \min(a,2n-a)\\i \equiv a \mod 2}} \Lo^{\frac{a-i}{2}} \pi_{i}
\end{equation}
with $\pi_i \in \largewedge^i W^*$ primitive. We define for $0 \leq i \leq \min(a,2n-a), i \equiv a \mod 2$ the maps 
\begin{displaymath}
 \Pi_i:\largewedge^aW^* \to \largewedge^i W^*, \phi \mapsto \pi_i.
\end{displaymath}

Let us fix a Euclidean scalar product on $W$.
Let 
\begin{displaymath}
 \Lambda:\largewedge^* W^* \to \largewedge^{*-2}W^*
\end{displaymath}
be the dual Lefschetz operator. Note that the two operators $\Lo$ and $\Lambda$ on $\largewedge^*W^*$, together with the counting operator $\mathrm{H}\phi=(a-n)\phi$ for $\phi \in \largewedge^a W^*$, define a representation of $\mathfrak{sl}_2$ on $\largewedge^*W^*$, in particular $[\Lo,\Lambda]=2\mathrm{H}$. An element $\pi \in \largewedge^iW^*$ is primitive if and only if $\Lambda \pi=0$. For more information we refer to \cite{huybrechts05}.

\begin{Lemma} \label{lemma_sl2_stuff}
\begin{enumerate}
 \item There are constants $c_{n,a,i}$ such that on $\largewedge^a W^*$
\begin{displaymath}
 \Lambda^s=\sum_{\substack{0 \leq i \leq \min(a-2s,2n-a)\\ i \equiv a \mod 2}}  c_{n,a,i} \Lo^{\frac{a-i}{2}-s} \circ \Pi_i.
\end{displaymath} 
\item For $i \leq a, i \equiv a \mod 2$ there are constants $c'_{n,a,i,j}$ such that on $\largewedge^a W^*$
\begin{displaymath}
 \Pi_i=\sum_{\max\left(0,\frac{a-i}{2}\right) \leq j \leq \frac{a}{2}} c'_{n,a,i,j} \Lo^{j+\frac{i-a}{2}} \circ \Lambda^j.
\end{displaymath}

\item Let $s \geq 0$. There are constants $c''_{n,a,s,j}$ such that on $\largewedge^a W^*$
\begin{displaymath}
 \Lambda^s = \sum_{j=\max(s,a-n)}^{\lfloor \frac{a}{2} \rfloor} c''_{n,a,s,j} \Lo^{j-s} \circ \Lambda^j.
\end{displaymath}
\end{enumerate}
\end{Lemma}

\proof
\begin{enumerate}
 \item The first statement follows by induction on $s$, using the relation 
\begin{equation} \label{eq_huyb_cor}
\Lambda \circ \Lo^j \pi_i=[\Lambda,\Lo^j]\pi_i=-j(i-n+j-1)\Lo^{j-1}\pi_i,
\end{equation}
which is a consequence of the commutator relations (compare also \cite[Corollary 1.2.28]{huybrechts05}).

\item We use induction on $i$, with the empty case $i<0$ the induction start. Suppose that we have shown the statement for all $i < i_0$, where $i_0 \equiv a \mod 2$.
Using \eqref{eq_lefschetz_decomp} and repeated application of \eqref{eq_huyb_cor}, we compute 
\begin{align*}
 \Lambda^\frac{a-i_0}{2} \phi & =\sum_{\substack{0 \leq i \leq \min(a,2n-a)\\i \equiv a \mod 2}} \Lambda^\frac{a-i_0}{2} \Lo^{\frac{a-i}{2}} \pi_{i}\\
 & = \sum_{\substack{0 \leq i \leq \min(a,2n-a,i_0)\\i \equiv a \mod 2}} \tilde c_{n,a,i,i_0} \Lo^\frac{i_0-i}{2}\pi_i,
\end{align*}
where $\tilde c_{n,a,i_0,i_0} \neq 0$. Using the induction hypothesis, we see that $\pi_{i_0}$ can be written as a combination of $\Lo^{j+\frac{i_0-a}{2}} \circ \Lambda^j \phi$, with coefficients independent of $\phi$.
\item The last equation is a direct consequence of the first two.

\end{enumerate}
\endproof

\begin{Proposition} \label{prop_fft} 
For $0 \leq a,b \leq 2n$ and $r \geq 0$, consider the space 
\begin{displaymath}
  \Xi:=\Hom_{G_0}(\largewedge^aW^* \otimes \Sym^rW^*,\largewedge^bW^*).
\end{displaymath}
\begin{enumerate}
\item If $r+a+b$ is odd, or if $r>2$, then $\Xi=0$.
 \item If $r=0$ and $a+b$ is even, then $\Xi$ is generated by all maps 
\begin{equation}   \label{eq_caser0}
\Lo^t \circ \Lambda^s, \quad \max(0,a-n) \leq s \leq \frac{a}{2}, 0 \leq t \leq \frac{b}{2}, t-s=\frac{b-a}{2}.
\end{equation}
\item If $r=1$ and $a \equiv b+1 \mod 2$, then $\Xi$ is generated by the maps 
\begin{equation} \label{eq_caser1_type1}
 \phi \otimes \tau \mapsto \Lo^t \circ \Lambda^s (\phi \wedge \tau),
\end{equation}
where $\max(0,a+1-n) \leq s \leq \frac{a+1}{2}, 0 \leq t \leq \frac{b}{2}, t-s=\frac{b-a-1}{2}$;
and the maps 
\begin{equation} \label{eq_caser1_type2}
 \phi \otimes \tau \mapsto \Lo^t \circ \Lambda^s (\phi) \wedge \tau, 
\end{equation}
where  $\max(0,a-n) \leq s \leq \frac{a}{2}, 0 \leq t \leq \frac{b-1}{2},t-s=\frac{b-a-1}{2}$.
\item If $r=2$ and $a \equiv b \mod 2$, then $\Xi$ is generated by all maps
\begin{displaymath}
\phi \otimes \tau_1 \otimes \tau_2 \mapsto \Lo^t \circ \Lambda^s(\phi \wedge \tau_1) \wedge \tau_2+ L^t \circ \Lambda^s(\phi \wedge \tau_2) \wedge \tau_1,
\end{displaymath}
where $\max(0,a+1-n) \leq s \leq \frac{a+1}{2}, 0 \leq t \leq \frac{b-1}{2},t-s=\frac{b-a-2}{2}$. 
\end{enumerate}
\end{Proposition}

\proof
Let $E,F$ be linear representations of $G_0$ and $E' \subset E$ a subrepresentation. Any $G_0$-equivariant linear map $\rho:E' \to F$ can be lifted to an equivariant map $\tilde \rho:E \to F$, as $G_0$ is semi-simple. 

Let $\iota:\largewedge^aW^* \otimes \Sym^rW^* \to (W^*)^{\otimes (a+r)}$ be the natural inclusion map and $\pi:(W^*)^{\otimes b} \to \largewedge^bW^*$ be the natural projection map. Then any equivariant map $\rho:\largewedge^aW^* \otimes \Sym^rW^* \to \largewedge^bW^*$ can be lifted to an equivariant map $\tilde \rho:(W^*)^{\otimes (a+r)} \to (W^*)^{\otimes b}$, i.e. the following diagram commutes
\begin{displaymath}
 \xymatrix{(W^*)^{\otimes (a+r)} \ar[r]^-{\tilde \rho} & (W^*)^{\otimes b} \ar[d]^\pi\\
 \largewedge^aW^* \otimes \Sym^rW^* \ar[r]^-\rho \ar[u]^\iota & \largewedge^bW^*
  }
\end{displaymath}

Note that $\Omega$ induces an equivariant identification $W \cong W^*$. Under this identification, $\Omega \in \largewedge^2W^* \subset W^* \otimes W^*$ corresponds to the element $\sum_i e_i \wedge f_i \in \largewedge^2 W \subset W \otimes W$, where $\{e_1,f_1,\ldots,e_n,f_n\}$ is a symplectic basis of $(W,\Omega)$ (i.e. $\Omega(e_i,f_j)=\delta_{ij}, \Omega(e_i,e_j),\Omega(f_i,f_j)=0$). There are contractions $W \otimes W \to \C, v \otimes w \mapsto \Omega(v,w)$ and $W^* \otimes W^* \to \C, \omega \otimes \tau \mapsto \sum_i [\omega(e_i)\tau(f_i)-\omega(f_i)\tau(e_i)]$. The operator $\Lambda:\largewedge^*W^* \to \largewedge^{*-2}W^*$ is the composition of the contraction followed by antisymmetrization.

Set $a':=a+r$ and $m:=a'+b$. The space $\Hom_{G_0}((W^*)^{\otimes a'},(W^*)^{\otimes b})$ may be identified with the space of invariant elements in $(W^*)^{\otimes m}$. It is trivial if $m$ is odd. If $m=2l$ is even, this space is spanned by the complete contractions, i.e. maps of the form 
\begin{displaymath}
 \rho(w_1,\ldots,w_m)=\Omega(w_{i_1},w_{j_1}) \cdots \Omega(w_{i_l},w_{j_l}),
\end{displaymath}
where $\{(i_1,j_1),\ldots,(i_l,j_l)\}$ is a $2$-partition of $\{1,\ldots,m\}$. This follows from the first fundamental theorem of invariant theory for the group $G_0$ (see \cite[Theorem 5.2.2]{goodman_wallach09}).  

Let us describe this map, unwinding the identifications between $W$ and $W^*$. Fix some $2$-partition as above. Let $s$ be the number of pairs $(i,j)$ such that both $i$ and $j$ are in $\{1,\ldots,a'\}$; let $t$ be the number of pairs such that both $i$ and $j$ are in $\{a'+1,\ldots,b\}$ and let $u$ be the number of pairs such that one entry belongs to $\{1,\ldots,a'\}$ and the other entry to $\{a'+1,\ldots,b\}$. Clearly $u+2s=a',u+2t=b$. The corresponding invariant element, considered as a map $(W^*)^{\otimes a'} \to (W^*)^{\otimes b}$ can be described as follows. First contract $(W^*)^{\otimes a'}$ in the places given by the $s$ pairs of the first type. We get a map $(W^*)^{\otimes a'} \to (W^*)^{\otimes (a'-2s)}$. Then identify $(W^*)^{\otimes (a'-2s)}$ and $(W^*)^{\otimes (b-2t)}$ by using the $u$ pairs of the third type. Then multiply by $\Omega$  in each of the $2t$ places given by the $t$ pairs of the second type.
We thus get a map
\begin{displaymath}
\tilde \rho: W^{\otimes a'} \to W^{\otimes (a'-2s)} \to W^{\otimes (b-2t)} \to W^{\otimes b}.
\end{displaymath}

The corresponding map $\rho:\largewedge^aW^* \otimes \Sym^r W^* \to \largewedge^b W^*$ is given by $\pi \circ \tilde \rho \circ \iota$. Suppose that there is a pair $(i,j)$ with both $i$ and $j$ in $\{a+1,\ldots,a+r\}$. Since $\pi$ is an antisymmetrization, while $\Sym^rW^*$ is symmetric, $\rho=0$ in this case. Similarly, suppose that there are two pairs $(i,j), (i',j')$ such that $i,i' \in \{1,\ldots,a\}, j,j' \in \{a+1,\ldots,a+r\}$. Then by symmetry and antisymmetry, $\rho=0$. This proves that all maps are trivial if $r>2$, since in this case there are necessarily pairs of one of these kinds. 

If $r=0$, the map $\rho$ is given by $\phi \mapsto \Lo^t \circ \Lambda^s \phi$. Hence $\Xi$ is in this case spanned by the maps $\Lo^t \circ \Lambda^s$ with $0 \leq s \leq \frac{a}{2}, 0 \leq t \leq \frac{b}{2}, t-s=\frac{b-a}{2}$. By Lemma \ref{lemma_sl2_stuff}, we may even further restrict to $s \geq a-n$.

If $r=1$, we either get a pair $(i,a+1)$ with $i \in \{1,\ldots,a\}$ or a pair $(i,a+1)$ with $i \in \{a'+1,\ldots,m\}$. The map $\rho$ is therefore of one of two kinds: $\phi \otimes \tau \mapsto \Lo^t \circ \Lambda^s (\phi \wedge \tau)$ with $0 \leq s \leq \frac{a+1}{2}, 0 \leq t \leq \frac{b}{2}, t-s=\frac{b-a-1}{2}$ or $\phi \otimes \tau \mapsto \Lo^t \circ \Lambda^s (\phi) \wedge \tau$ with $0 \leq s \leq \frac{a}{2}, 0 \leq t \leq \frac{b-1}{2}, t-s=\frac{b-a-1}{2}$. Lemma \ref{lemma_sl2_stuff} allows to assume moreover $s \geq a+1-n$ in the first case and $s \geq a-n$ in the second case.

If $r=2$, we may assume that we have one pair $(i,a+1)$ with $i \in \{1,\ldots,a\}$ and one pair $(j,a+2)$ with $j \in \{a'+1,\ldots,m\}$. Then the map $\rho$ is given by $\phi \otimes \tau_1 \otimes \tau_2 \mapsto \Lo^t \circ \Lambda^s(\phi \wedge \tau_1) \wedge \tau_2+\Lo^t \circ \Lambda^s(\phi \wedge \tau_2) \wedge \tau_1$ with $0 \leq s \leq \frac{a+1}{2}, 0 \leq t \leq \frac{b-1}{2},t-s=\frac{b-a-2}{2}$. By Lemma \ref{lemma_sl2_stuff}, we may even further restrict to $s \geq a+1-n$.
\endproof

Sometimes it will be easier to work with the projections $\Pi_i$ instead of the iterated dual Lefschetz operator $\Lambda^s$. We therefore reformulate the preceding proposition in terms of the $\Pi_i$'s. 

\begin{Corollary} \label{cor_pi_instead_of_lambda}
\begin{enumerate}
\item If $r=0$ and $a+b$ is even, then $\Xi$ is generated by all maps 
\begin{displaymath} 
 \Lo^{\frac{b-i}{2}} \circ \Pi_i,
\end{displaymath}
where $0 \leq i \leq \min(a,2n-a,b,2n-b), i \equiv b \mod 2$.
\item If $r=1$ and $a \equiv b+1 \mod 2$, then $\Xi$ is generated by the maps 
\begin{displaymath} 
 \phi \otimes \tau \mapsto \Lo^{\frac{b-i}{2}} \circ \Pi_i(\phi \wedge \tau),
\end{displaymath}
with $0 \leq i \leq \min(a+1,2n-a-1,b,2n-b), i \equiv b \mod 2$ and the maps 
\begin{displaymath}
 \phi \otimes \tau \mapsto \Lo^{\frac{b-1-i}{2}} \circ \Pi_i(\phi) \wedge \tau,
\end{displaymath}
with $0 \leq i \leq \min(a,2n-a,b-1,2n-b+1), i \equiv b+1 \mod 2$.
\item If $r=2$ and $a \equiv b \mod 2$, then $\Xi$ is generated by all maps
\begin{displaymath}
 \phi \otimes \tau_1 \otimes \tau_2 \mapsto \Lo^{\frac{b-i-1}{2}} \left[\Pi_i(\phi \wedge \tau_1) \wedge \tau_2 + \Pi_i(\phi \wedge \tau_2) \wedge \tau_1\right],
\end{displaymath}
where $0 \leq i \leq \min(a+1,2n-a-1,b-1,2n-b+1), i \equiv b+1 \mod 2$. 
\end{enumerate}
\end{Corollary}

\proof
Follows from Proposition \ref{prop_fft} by applying Lemma \ref{lemma_sl2_stuff} (i).
\endproof

\section{Reduction to the horizontal part}
\label{sec_reduction_horizontal}

Recall that $V=T_0\R^{2n+1}=W \oplus D$ with $W=Q_0 \cong \C^n, D \cong \R$. Moreover, we have seen in Section \ref{sec_isotropy_group} that the additive group $W=\C^n$ is a subgroup of $G$. Correspondingly, we obtain an action by the Lie algebra $W$ (with trivial brackets) on $V$. For $w \in W$, we write $R_w \in \mathfrak{g}$ for the corresponding element. We also get induced actions on the spaces $\largewedge^aV^*$ and $\Sym^rV^*$. 

\begin{Proposition} \label{prop_reduction_horizontal}
\begin{enumerate}
\item Let $a<2n$. Let $E \subset \largewedge^aV^* \otimes \Sym^ rV^ *$ be a subspace such that $R_w \phi \in E$ for each $w \in W, \phi \in E$. If $E$ contains $\largewedge^aW^* \otimes \Sym^ rW^ *$, then $E=\largewedge^aV^* \otimes \Sym^ rV^ *$.
\item Let $E \subset \largewedge^{2n}V^* \otimes \Sym^ rV^ *$ be a subspace such that $R_w \phi \in E$ for each $w \in W, \phi \in E$. If $E$ contains $\largewedge^{2n}W^* \otimes \Sym^ rV^ *$, then $E=\largewedge^{2n}V^* \otimes \Sym^ rV^ *$.
\end{enumerate}
\end{Proposition}

\proof
\begin{enumerate}
\item Let $\{e_1,f_1,\ldots,e_n,f_n\}$ be a symplectic basis of $W$ and $0 \neq h \in D$. Then $\{e_1,f_1,\ldots,e_n,f_n,h\}$ is a basis of $V$ and we denote the dual basis by $\{e_1^*,f_1^*,\ldots,e_n^*,f_n^*,h^*\}$. If $\xi \in V^*$ and $w \in W$, then $R_w(\xi)=\xi(w)h^*$ (strictly speaking, this equality holds up to a scalar, but we may scale $h$ in such a way that the scalar equals $1$). 

We decompose 
\begin{displaymath}
\largewedge^aV^* \otimes \Sym^ rV^ *= \bigoplus_{\epsilon=0,1} \bigoplus_{l=0}^r \largewedge^{a-\epsilon}W^* \otimes (D^*)^\epsilon \otimes \Sym^{r-l}W^* \otimes (D^*)^l,
\end{displaymath}
and show that each summand belongs to $E$.

Take $I,J \subset \{1,\ldots,n\}$ with $\#I+\#J=a<2n$. Suppose without loss of generality that $1 \notin I$. For $0 \leq l \leq r$, let $\tau_1,\ldots,\tau_{r-l} \in \{e_1^*,\ldots,e_n^*,f_1^*,\ldots,f_n^*\}$ and let $z \geq 0$ be the number of times $e_1^*$ appears among these elements. 
Then 
\begin{displaymath}
 R_{e_1}^l(e_I^* \wedge f_J^* \otimes \tau_1 \cdot \ldots \cdot \tau_{r-l} \cdot (e_1^*)^l)=\frac{(l+z)!}{z!} e_I^* \wedge f_J^* \otimes \tau_1 \cdot \ldots \cdot \tau_{r-l} \otimes (h^*)^l \in E.
\end{displaymath}
This shows that all summands with $\epsilon=0$ belong to $E$.

Next, let $e_I^* \wedge f_J^* \in \largewedge^{a-1}W^*$ and $\tau_1,\ldots,\tau_r \in V^*$. Suppose again that $1 \notin I$. Then $e_1^* \wedge e_I^* \wedge f_J^* \otimes \tau_1 \cdot \ldots \cdot \tau_r \in E$ and hence 
\begin{equation}
 R_{e_1}(e_1^* \wedge e_I^* \wedge f_J^* \otimes \tau_1 \cdot \ldots \cdot \tau_r)=e_I^* \wedge f_J^* \otimes h^* \otimes \tau_1 \cdot \ldots \cdot \tau_r+e_1^* \wedge e_I^* \wedge f_J^* \otimes R_{e_1}(\tau_1 \cdot \ldots \cdot \tau_r).
\end{equation}
Since $e_1^* \wedge e_I^* \wedge f_J^* \otimes R_{e_1}(\tau_1 \cdot \ldots \cdot \tau_r) \in \largewedge^aW^* \otimes \Sym^rV^*$, this term belongs to $E$ by what we have already shown. 
It follows that $e_I^* \wedge f_J^* \otimes h^* \otimes \tau_1 \cdot \ldots \cdot \tau_r \in E$.
\item We argue as in the second part of the previous proof. This time the terms $e_1^* \wedge e_I^* \wedge f_J^* \otimes \tau_1 \cdot \ldots \cdot \tau_r$ and $e_1^* \wedge e_I^* \wedge f_J^* \otimes R_{e_1}(\tau_1 \cdot \ldots \cdot \tau_r)$ belong to $E$ by assumption.
\end{enumerate}
\endproof

\begin{Corollary} \label{cor_reduction_horizontal}
Let
\begin{displaymath}
 \sigma:\largewedge^a V^* \otimes \Sym^r V^* \to \largewedge^bV^*
\end{displaymath}
be a $G$-equivariant linear map. 
\begin{enumerate}
 \item If $a<2n$ and if $\sigma$ vanishes on $\largewedge^aW^* \otimes \Sym^rW^*$, then $\sigma=0$. 
 \item If $a=2n$ and if $\sigma$ vanishes on $\largewedge^aW^* \otimes \Sym^rV^*$, then $\sigma=0$.
\end{enumerate}
\end{Corollary}

\proof
Since $\sigma$ commutes with each $R_w$, its kernel is invariant under each $R_w$. The statement therefore follows from Proposition \ref{prop_reduction_horizontal}.
\endproof

\section{Operators of degree $\leq 2$}
\label{sec_operators_degree2}

In this section we prove Proposition \ref{prop_invariant_symbols}, which will complete the proof of Theorem \ref{mainthm_fixed_dimension}. 

\proof[Proof of Proposition \ref{prop_invariant_symbols}]
Let $V:=T_0\R^{2n+1}, W:=Q_0$. Let $\sigma:\largewedge^aV^* \otimes \Sym^r V^*\to \largewedge^bV^*$ be $G$-equivariant. We have to show that  $\sigma_{\Po,0}=\sigma$, where $\Po \in R_{a,b}$.

Let us first assume that $a<2n$. We can decompose $\largewedge^bV^*=\largewedge^bW^* \oplus \largewedge^{b-1} \otimes D^*$ and look at the corresponding two components of $\sigma$. By homogeneity, the restriction to $\largewedge^aW^* \otimes \Sym^r W^*$ vanishes if $b \neq a+r,a+r-1$, and Corollary \ref{cor_reduction_horizontal} implies that $\sigma=0$ in this case. 

If $b=a+r$, we obtain by restriction a $G_0$-equivariant map $\sigma:\largewedge^aW^* \otimes \Sym^r W^* \to \largewedge^{a+r}W^*$, and Proposition \ref{prop_fft} classifies these maps. Similarly, if $b=a+r-1$, we obtain a $G_0$-equivariant map $\sigma:\largewedge^aW^* \otimes \Sym^r W^* \to \largewedge^{a+r-2}W^* \otimes D^* \cong \largewedge^{a+r-2} W^*$. In both case, $\sigma$ is trivial if $r>2$. 

However, even for $r \leq 2$, not all $G_0$-equivariant maps from $\largewedge^aW^* \otimes \Sym^r W^*$ to $\largewedge^{a+r}W^*$ or $\largewedge^{a+r-2} W^*$ are restrictions of $G$-equivariant maps $\largewedge^aV^* \otimes \Sym^r V^*$ to $\largewedge^bV^*$. To rule out those that are not, we use the easy observation that if $\sigma$ is $G$-equivariant, then, for each fixed $w \in W$, $\sigma$ must map the kernel of $R_w$ to itself, since $R_w \circ \sigma=\sigma \circ R_w$. 

\begin{enumerate}
 \item {\bf Case $r=0, b=a-1$.}
 The restriction of $\sigma$ to the horizontal part $W$ is an element of 
 \begin{displaymath}
\Hom_{G_0}(\largewedge^aW^*,\largewedge^{a-2}W^* \otimes D^*) \cong \Hom_{G_0}(\largewedge^aW^*,\largewedge^{a-2}W^*). 
 \end{displaymath} 
 By Corollary \ref{cor_pi_instead_of_lambda}, $\sigma|_W$ is a linear combination of the operators $\Lo^{\frac{a-i-2}{2}} \circ \Pi_i \wedge \alpha$, where $0 \leq i \leq \min(2n-a,a-2), i \equiv a \mod 2$, i.e. a linear combination of the $\po_{a,i}$. Since $\po_{a,i}$ is the symbol 
 of $\Po_{a,i}$, we are done in this case.
 
\item {\bf Case $r=0, b=a$.}
Any $G_0$-invariant operator $\sigma \in \Hom_{G_0}(\largewedge^aW^*,\largewedge^{a}W^*)$ can be written as 
\begin{displaymath}
\sigma(\phi)=\sum_{j=\max(0,a-n)}^{\left\lfloor \frac{a}{2}\right\rfloor} c_j \Lo^j \Lambda^j(\phi), \quad \phi \in \largewedge^a W^*,
\end{displaymath}
for some constants $c_j$. 

The identity $\largewedge^aW^* \to \largewedge^aW^*$ can be expressed as some linear combination: $\phi=\sum_{j=\max(0,a-n)}^{\left\lfloor \frac{a}{2}\right\rfloor} \tilde c_j \Lo^j \Lambda^j(\phi)$ for all $\phi \in \largewedge^a W^*$, with $\tilde c_{\max(0,a-n)} \neq 0$. Since the identity is equivariant, we may subtract some multiple of it from $\sigma$, and hence assume that $c_{\max(0,a-n)}=0$. We now prove by induction that $c_j=0$ for all $j>\max(0,a-n)$. Suppose that we have shown this equality for all indices strictly smaller than $j$. 

Take 
\begin{displaymath}
 \phi:=(e_2^* \wedge f_2^*) \wedge \ldots \wedge (e_{j+1}^* \wedge f_{j+1}^*) \wedge e_{j+2}^* \wedge \ldots \wedge e_{a-j+1}^* \in \largewedge^aW^*.
\end{displaymath}

Then $\Lambda^j\phi= j! e_{j+2}^* \wedge \ldots \wedge e_{a-j+1}^* \neq 0$ and $\Lambda^{j+1}\phi=0$. By induction hypothesis we obtain that 
\begin{displaymath}
\sigma \phi=c_j \Lo^j \Lambda^j \phi,
\end{displaymath}
and hence 
\begin{align*}
 0=\sigma R_{e_1} \phi=R_{e_1} \sigma \phi & = j c_j h^* \wedge f_1^* \wedge \Lo^{j-1} \circ \Lambda^j \phi.
\end{align*}
It follows that $c_j=0$, as claimed.  

\item {\bf Case $r=1, b=a+1$.} The restriction $\sigma \in \Hom_{G_0}(\largewedge^aW^* \otimes W^*,\largewedge^{a+1} W^*)$ is given by
\begin{displaymath}
 \phi \otimes \tau \mapsto \sum_{j=\max(0,a-n+1)}^{\left\lfloor \frac{a+1}{2}\right\rfloor} c_j \Lo^j \Lambda^j(\phi \wedge \tau)+\sum_{j=\max(0,a-n)}^{\left\lfloor \frac{a}{2}\right\rfloor}  \tilde c_j \Lo^j \Lambda^j(\phi) \wedge \tau
\end{displaymath}
for some constants $c_j,\tilde c_j$. The map $\phi \otimes \tau \mapsto \phi \wedge \tau$ may be written as 
\begin{align*}
 \phi \wedge \tau & = \sum_{j=\max(0,a-n+1)}^{\left\lfloor \frac{a+1}{2}\right\rfloor} c_j' \Lo^j \Lambda^j(\phi \wedge \tau)\\
 \phi \wedge \tau & = \sum_{j=\max(0,a-n)}^{\left\lfloor \frac{a}{2}\right\rfloor}  \tilde c_j' \Lo^j \Lambda^j(\phi) \wedge \tau,
\end{align*}
where $c'_{\max(0,a-n+1)} \neq 0, \tilde c'_{\max(0,a-n)} \neq 0$. Subtracting such terms from $\sigma$, we may assume that $c_{\max(0,a-n+1)} = 0, \tilde c_{\max(0,a-n)} = 0$. We now prove by induction on $j$ that $c_j=0$ for all $j>\max(0,a-n+1)$ and $\tilde c_j=0$ for all $j>\max(0,a-n)$. Suppose that we have shown these equalities for all indices strictly smaller than $j$. 

Let $j>a-n$ and take 
\begin{align*}
 \phi & :=(e_2^* \wedge f_2^*) \wedge \ldots \wedge (e_{j+1}^* \wedge f_{j+1}^*) \wedge e_{j+2}^* \wedge \ldots \wedge e_{a-j+1}^* \in \largewedge^aW^*\\
 \tau & := e_{a-j+1}^* \in W^*,
\end{align*}
which is possible since $a-j+1 < a-(a-n)+1=n+1$.

Then $\sigma(\phi \otimes \tau)=\tilde c_j (\Lo^j \Lambda^j\phi) \wedge \tau$ and hence 
\begin{align*}
0=\sigma R_{e_1}(\phi \otimes \tau)= R_{e_1} \sigma (\phi \otimes \tau)& = j \tilde c_j h^* \wedge f_1^* \wedge \Lo^{j-1} \circ \Lambda^j \phi \wedge e_{a-j+1}^*.
\end{align*}
It follows that $\tilde c_j=0$.  

To obtain a second equation, let $j \geq a-n+2$ (hence $a-j+2 \leq n$) and put
\begin{align*}
 \phi & :=(e_2^* \wedge f_2^*) \wedge \ldots \wedge (e_j^* \wedge f_j^*) \wedge e_{j+1}^* \wedge \ldots \wedge e_{a-j+2}^* \in \largewedge^aW^*\\
 \tau & := f_{j+1}^* \in W^*.
\end{align*}
Then $\Lambda^j\phi=0$ but $\Lambda^j(\phi \wedge \tau) \neq 0$. We obtain $\sigma(\phi \otimes \tau)=c_j \Lo^j \Lambda^j(\phi \wedge \tau)$ and 
\begin{displaymath}
 0=\sigma R_w(\phi \otimes \tau)=R_w \sigma (\phi \otimes \tau) = j c_j h^* \wedge f_1^* \wedge \Lo^{j-1} \circ \Lambda^j (\phi \wedge \tau).
\end{displaymath}
It follows that $c_j=0$, which finishes the induction. 

We conclude that $\sigma$ is a multiple of the map $\phi \otimes \tau \mapsto \phi \wedge \tau$, which is the (horizontal part of) the symbol of the usual exterior derivative $\ed$.  

\item {\bf Case $r=1,b=a$.} By Corollary \ref{cor_pi_instead_of_lambda}, $\sigma \in \Hom_{G_0}(\largewedge^a W^* \otimes W^*,\largewedge^{a-1} W^* \otimes D^*) \cong \Hom_{G_0}(\largewedge^a W^* \otimes W^*,\largewedge^{a-1} W^*)$ is a linear combination of the maps   
 \begin{displaymath}
 \phi \otimes \tau \mapsto \Lo^\frac{a-i-1}{2} \circ \Pi_i(\phi \wedge \tau) \wedge \alpha, 0 \leq i \leq \min(a-1,2n-a-1), i \equiv a+1 \mod 2.
 \end{displaymath}
 and 
 \begin{displaymath}
 \phi \otimes \tau \mapsto \Lo^\frac{a-i-2}{2} \circ \Pi_i(\phi) \wedge \tau \wedge \alpha, 0 \leq i \leq \min(a-2,2n-a), i \equiv a \mod 2. 
 \end{displaymath}
They equal the horizontal part of the symbols of $\Po_{a+1,i} \circ \ed$ and $\ed \circ \Po_{a,i}$. 
\item {\bf Case $r=2,b=a+2$.} 
By Proposition \ref{prop_fft}, the restriction $\sigma \in \Hom_{G_0}(\largewedge^aW^* \otimes \Sym^2W^*,\largewedge^{a+2}W^*)$ is given by
\begin{displaymath}
\phi \otimes \tau_1 \otimes \tau_2 \mapsto \sum_{j=\max(1,a-n+2)}^{\left\lfloor \frac{a+1}{2}\right\rfloor} c_j \left(\Lo^j \Lambda^j(\phi \wedge \tau_1) \wedge \tau_2+ \Lo^j \Lambda^j(\phi \wedge \tau_2) \wedge \tau_1\right).
\end{displaymath}

We prove by induction on $j$ that $c_j=0$.
 
Take 
\begin{align*}
\phi & :=(e_2^* \wedge f_2^*) \wedge \ldots \wedge (e_j^* \wedge f_j^*) \wedge e_{j+1}^* \wedge \ldots \wedge e_{a-j+2}^* \in \largewedge^aW^*\\
\tau_1 & := e_{j+1}^* \in W^*,\\
\tau_2 & := f_{j+1}^* \in W^*,
\end{align*}
which is possible since $a-j+2 \leq a-(a-n+2)+2=n$.

Then 
\begin{align*}
\Lambda^j(\phi \wedge \tau_1) & = 0 \\
\Lo^j\circ \Lambda^j(\phi \wedge \tau_2) \wedge \tau_1 & = \pm j! \Lo^j(e_{j+1}^* \wedge \ldots \wedge e_{a-j+2}^*) 
\end{align*}
and hence 
\begin{displaymath}
0=\sigma \circ R_{e_1} (\phi \otimes \tau_1 \otimes \tau_2)= R_{e_1} \circ \sigma (\phi \otimes \tau_1 \otimes \tau_2)=\pm j! c_j  R_{e_1} \Lo^j(e_{j+1}^* \wedge \ldots \wedge e_{a-j+2}^*). 
\end{displaymath}
It follows that $c_j=0$.  
\item {\bf Case $r=2,b=a+1$.} By Corollary \ref{cor_pi_instead_of_lambda}, the restriction $\sigma \in \Hom_{G_0}(\largewedge^aW^* \otimes \Sym^2W^*,\largewedge^a W^* \otimes D^*) \cong \Hom_{G_0}(\largewedge^aW^* \otimes \Sym^2W^*,\largewedge^a W^*)$ is given by a linear combination of the maps 
\begin{align*}
 \phi \otimes \tau_1 \otimes \tau_2 & \mapsto [\Lo^\frac{a-i-1}{2} \Pi_i(\phi \wedge \tau_1) \wedge \tau_2+\Lo^\frac{a-i-1}{2} \Pi_i(\phi \wedge \tau_2) \wedge \tau_1], 
\end{align*}
where $0 \leq i \leq \min(a-1,2n-a-1), i \equiv a+1 \mod 2$. This is the horizontal part of the symbol of $\ed \circ \Po_{a+1,i} \circ \ed$.
\end{enumerate}

To complete the proof, we have to study the cases $a=2n,2n+1$.

Let $a=2n$. 

By Corollary \ref{cor_reduction_horizontal}ii), it is enough to prove that if $\sigma:\largewedge^{2n}W^* \otimes \Sym^rV^* \to \largewedge^bV^*$ is $G$-equivariant, then there exists $\Po \in R_{2n,b}^{2n+1}$ with $\sigma_{\Po,0}=\sigma$. Recall that $R_{2n,b}^{2n+1}$ is spanned by $\Po_{2n,0}$ if $b=2n-1$; by $\ed \circ \Po_{2n,0}$ and $\mathrm{id}$ if $b=2n$; and by $\ed$ if $b=2n+1$. 

We decompose 
\begin{displaymath}
 \largewedge^{2n}W^* \otimes \Sym^rV^*=\bigoplus_{l=0}^r \largewedge^{2n}W^* \otimes \Sym^{r-l}W^* \otimes (D^*)^{\otimes l}
\end{displaymath}
and 
\begin{displaymath}
 \largewedge^bV^*=\bigoplus_{\epsilon=0}^1 \largewedge^{b-\epsilon}W^* \otimes (D^*)^{\otimes \epsilon},
\end{displaymath}
and decompose $\sigma$ according to these splittings. By looking at degrees of homogeneity, we must have $2n+r+l=b+\epsilon$. Since $b \leq 2n+1$, we are left with the following cases
\begin{enumerate}
 \item $r=2,l=0,b=2n+1,\epsilon=1$. Then $\sigma:\largewedge^{2n}W^* \otimes \Sym^2 W^* \to \largewedge^{2n}W^* \otimes D^*$ is zero by Proposition \ref{prop_fft} iv). 
 \item $r=1,l=1,b=2n+1,\epsilon=1$. Then $\sigma:\largewedge^{2n}W^* \otimes D^* \to \largewedge^{2n}W^* \otimes D^*$ is a multiple of the identity by Corollary \ref{cor_pi_instead_of_lambda}i), and the identity is the restriction of the symbol of $\ed$ to $\largewedge^{2n}W^* \otimes D^*$.
 \item $r=1,l=0,b=2n,\epsilon=1$. Then $\sigma:\largewedge^{2n}W^* \otimes W^* \to \largewedge^{2n-1}W^* \otimes D^*$ is a multiple of $\sigma_{\ed \circ \Po_{2n,0},p}$ by Corollary \ref{cor_pi_instead_of_lambda}ii).
 \item $r=0,l=0,b=2n,\epsilon=0$. Then $\sigma:\largewedge^{2n}W^* \to \largewedge^{2n}W^*$ is a multiple of the identity, which is the restriction of the symbol of the identity map to $\largewedge^{2n}W^*$.  
 \item $r=0,l=0,b=2n-1,\epsilon=1$. By Corollary \ref{cor_pi_instead_of_lambda}i), $\sigma:\largewedge^{2n}W^* \to \largewedge^{2n-2}W^* \otimes D^*$ is the restriction of the symbol of $\Po_{2n,0}$ to $\largewedge^{2n}W^*$. 
\end{enumerate}

If $a=2n+1$, then by looking at degrees of homogeneity we find that $\sigma=0$ except if $r=0,b=2n+1$, in which case $\sigma$ is a multiple of the identity. 
\endproof

\def\cprime{$'$}

%\bibliographystyle{plain}
%\bibliography{../biblio}
\end{document}